# Système dérivé et suite duale d'une suite barypolygonale - Partie 2

Par David Pouvreau[1] et Vincent Bouis[2]

**Résumé**

Dans le prolongement de la première partie de cet article, les propriétés du système dérivé et de la suite duale d'une suite barypolygonale sont ici établies dans toute leur généralité : il est montré que les résultats obtenus dans le cas où $p = 3$ sont conservés si $p \geq 4$.

**Abstract**

Continuing the first part of this paper, the present one establishes in their full generality the properties of the derived system and of the dual sequence of any barypolygonal sequence: it is proven that the results obtained in the case $p = 3$ remain valid if $p \geq 4$.

## 1. Introduction

Cet article poursuit l'étude initiée dans (Pouvreau, Bouis, 2018), où les notions de système dérivé et de suite duale d'une suite barypolygonale ont été définies. Rappelons brièvement ces notions.

Soit $E$ un espace affine réel de dimension finie quelconque. Soit $p \in \mathbb{N}\backslash\{0; 1\}$. On considère une famille ordonnée $\mathcal{A}$ de points distincts $(A_k)_{1\leq k\leq p}$ de $E$ et une famille ordonnée $t = (t_k)_{1\leq k\leq p}$ de réels de $]0; 1[$. On note $\mathfrak{B}$ la suite $t$-barypolygonale de $\mathcal{A}$, définie et étudiée dans (Pouvreau, Eupherte, 2016). Soit $\left(t^{(m)}\right)_{m\in\mathbb{N}}$ la suite de $(]0; 1[^p)^{\mathbb{N}}$ définie par :

$$(S) : \begin{cases} t^{(0)} = t \\ \forall\, m \in \mathbb{N},\ t^{(m+1)} = \left(t_k^{(m+1)}\right)_{1\leq k\leq p} \quad \text{où } \forall\, k \in [\![1; p]\!],\ t_k^{(m+1)} = \prod_{\substack{1\leq i\leq p \\ i\neq k}} \left(1 - t_i^{(m)}\right) \end{cases}$$

$(S)$ est le « système barypolygonal dérivé » de $\mathfrak{B}$. Pour tout $m \in \mathbb{N}$, on peut alors définir au moyen de $t^{(m)}$ la « suite dérivée $m$-ième » de $\mathfrak{B}$ , notée $\mathfrak{B}^{(m)}$, comme la suite $t^{(m)}$–barypolygonale de $\mathcal{A}$.

Pour chaque $m \in \mathbb{N}$, la suite $\mathfrak{B}^{(m)}$ a un point limite que l'on note $G_m$, déterminé par :

$$G_m = \text{bar}\left\{\left(A_k; \prod_{\substack{1\leq i\leq p \\ i\neq k}} (1 - t_i^{(m)})\right)\right\}_{1\leq k\leq p}$$

La suite $(G_m)_{m\in\mathbb{N}}$ est appelée la suite duale de $\mathfrak{B}$.

La première partie de cet article a établi les propriétés du système dérivé et de la suite duale lorsque $p = 2$, lorsque $\mathfrak{B}$ est régulière (pour $p \geq 3$ quelconque) et dans tous les cas où $p = 3$ (que $\mathfrak{B}$ soit régulière ou non). Il s'agira ici de montrer que toutes les propriétés établies dans le cas où $p = 3$ se généralisent pour $p \geq 4$ ; en particulier, que la suite duale de $\mathfrak{B}$ converge vers le centre de gravité de $\mathcal{A}$.

[1] Professeur agrégé de mathématiques et docteur en histoire des sciences. Centre Universitaire de Mayotte et Institut Alexander Grothendieck (Université de Montpellier). Email : david_pouvreau@orange.fr

[2] Etudiant à l'Ecole Normale Supérieure de Paris. Email : vbouis@clipper.ens.fr

## 2. Dynamique des systèmes barypolygonaux dérivés et convergence duale pour tout $p \geq 3$

Une suite $t$–barypolygonale irrégulière $\mathfrak{B}$ de $\mathcal{A}$ étant donnée dans le cas où $\mathrm{Card}(\mathcal{A}) = p \geq 3$, on pose : $\forall\, k \in [\![1;p]\!],\ u_k^{(m)} = 1 - t_k^{(m)}$. Les $p$ suites $\left(u_k^{(m)}\right)_{m\in\mathbb{N}}$ sont les solutions du système récurrent « conjugué » du système dérivé $(S)$ de $\mathfrak{B}$ :

$$(\Sigma) : \begin{cases} \forall\, k \in [\![1;p]\!],\ u_k^{(0)} = 1 - t_k \\ \forall\, k \in [\![1;p]\!],\ \forall\, m \in \mathbb{N},\ u_k^{(m+1)} = 1 - \prod\limits_{\substack{1\leq i\leq p \\ i\neq k}} u_i^{(m)} \end{cases}$$

### 2.1. Unicité et instabilité exponentielle d'un point stationnaire

Commençons par établir le théorème suivant, qui vaut aussi bien pour les suites régulières :

**Théorème 1**
Soit $\mathfrak{B}$ une suite barypolygonale de $\mathcal{A}$, avec $\mathrm{Card}(\mathcal{A}) = p \geq 3$.
Le système dérivé de $\mathfrak{B}$ a un unique point stationnaire dans $]0;1[^p$ : celui dont toutes les coordonnées sont égales à $(1-\alpha_p)$, où $\alpha_p$ est l'unique solution dans $[0;1]$ de l'équation $(E_p)$ : $x^{p-1} + x - 1 = 0$. Ce point est exponentiellement instable.

<u>Démonstration</u>.
Notons $\Omega = (\omega_k)_{1\leq k\leq p} \in\, ]0;1[^p$. Supposons que $\Omega$ est un point stationnaire de $(\Sigma)$. Alors :

$$\forall\, k \in [\![1;p]\!],\ \omega_k = 1 - \prod_{\substack{1\leq i\leq p \\ i\neq k}} \omega_i$$

$$\omega_1 = 1 - \prod_{2\leq i\leq p} \omega_i \quad \text{et} \quad \forall\, k \in [\![2;p]\!],\ \omega_1 + \prod_{2\leq i\leq p} \omega_i = \omega_k + \prod_{\substack{1\leq i\leq p \\ i\neq k}} \omega_i$$

$$\Leftrightarrow \omega_1 = 1 - \prod_{2\leq i\leq p} \omega_i \quad \text{et} \quad \forall\, k \in [\![2;p]\!],\ (\omega_1 - \omega_k)\left(1 - \prod_{\substack{2\leq i\leq p \\ i\neq k}} \omega_i\right) = 0$$

$$\Leftrightarrow \omega_1 = 1 - \prod_{2\leq i\leq p} \omega_i \quad \text{et} \quad \forall\, k \in [\![2;p]\!],\ \omega_1 = \omega_k$$

On en déduit que chacun des $p$ réels $(\omega_k)_{1\leq k\leq p}$ est l'unique solution dans $]0;1[$ de l'équation $(E_p)$. Par conséquent : $\forall\, k \in [\![1;p]\!],\ \omega_k = \alpha_p$. Réciproquement, il est clair que le point $\Omega$ dont toutes les coordonnées sont égales à $\alpha_p$ est bien un point stationnaire de $(\Sigma)$.

Étudions maintenant la stabilité de ce point stationnaire. Notons à cette fin :

$$\forall\, k \in [\![1;p]\!],\ \forall\, m \in \mathbb{N},\ v_k^{(m)} = u_k^{(m)} - \alpha_p$$

On a alors :

$$\forall\, k \in [\![1;p]\!],\ \forall\, m \in \mathbb{N},\ v_k^{(m+1)} + \alpha_p = 1 - \prod_{\substack{1\leq i\leq p \\ i\neq k}} \left(v_i^{(m)} + \alpha_p\right)$$

Le système récurrent linéarisé au voisinage de l'origine s'écrit donc :

$$\forall\, k \in [\![1;p]\!],\ \forall\, m \in \mathbb{N},\ v_k^{(m+1)} = 1 - \alpha_p - {\alpha_p}^{p-1} - \sum_{\substack{1\le i\le p\\ i\ne k}} {\alpha_p}^{p-2}\, v_i^{(m)}$$

C'est-à-dire, compte tenu de la définition de $\alpha_p$ :

$$\forall\, k \in [\![1;p]\!],\ \forall\, m \in \mathbb{N},\ v_k^{(m+1)} = - \sum_{\substack{1\le i\le p\\ i\ne k}} {\alpha_p}^{p-2}\, v_i^{(m)}$$

La matrice de ce système linéaire est symétrique réelle, donc diagonalisable dans $\mathbb{R}$. Son polynôme caractéristique s'écrit, en notant $\beta_p = {\alpha_p}^{p-2}$ :

$$\chi_{M_p} = \begin{vmatrix} X & \beta_p & \beta_p & \cdots & \beta_p & \beta_p \\ \beta_p & X & \ddots & & & \beta_p \\ \beta_p & \beta_p & \ddots & \ddots & [\beta_p] & \vdots \\ \vdots & & \ddots & \ddots & \ddots & \vdots \\ \vdots & & [\beta_p] & \ddots & \ddots & \beta_p \\ \beta_p & \beta_p & \cdots & \cdots & \beta_p & X \end{vmatrix}_{[p]} = \begin{vmatrix} X-\beta_p & \beta_p & \beta_p & \cdots & \beta_p & \beta_p \\ \beta_p - X & X & \ddots & & & \beta_p \\ 0 & \beta_p & \ddots & \ddots & [\beta_p] & \vdots \\ \vdots & & \ddots & \ddots & \ddots & \vdots \\ \vdots & & [\beta_p] & \ddots & \ddots & \beta_p \\ 0 & \beta_p & \cdots & \cdots & \beta_p & X \end{vmatrix}_{[p]}$$

$$= (X-\beta_p) \begin{vmatrix} 1 & \beta_p & \beta_p & \cdots & \beta_p & \beta_p \\ -1 & X & \ddots & & & \beta_p \\ 0 & \beta_p & \ddots & \ddots & [\beta_p] & \vdots \\ \vdots & & \ddots & \ddots & \ddots & \vdots \\ \vdots & & [\beta_p] & \ddots & \ddots & \beta_p \\ 0 & \beta_p & \cdots & \cdots & \beta_p & X \end{vmatrix}_{[p]}$$

$$= (X-\beta_p) \begin{vmatrix} 1 & \beta_p & \beta_p & \cdots & \beta_p & \beta_p \\ 0 & X+\beta_p & 2\beta_p & \cdots & 2\beta_p & 2\beta_p \\ 0 & \beta_p & X & \beta_p & \cdots & \beta_p \\ \vdots & & \ddots & & \ddots\ [\beta_p] & \vdots \\ \vdots & & [\beta_p] & & \ddots\ X & \beta_p \\ 0 & \beta_p & \cdots & & \beta_p\ \ \beta_p & X \end{vmatrix}_{[p]}$$

$$= (X-\beta_p) \begin{vmatrix} X+\beta_p & 2\beta_p & 2\beta_p & \cdots & 2\beta_p \\ \beta_p & X & \beta_p & \cdots & \beta_p \\ \beta_p & \beta_p & \ddots & [\beta_p] & \vdots \\ \vdots & \ddots & [\beta_p] & \ddots & \beta_p \\ \beta_p & \beta_p & & \beta_p & X \end{vmatrix}_{[p-1]}$$

Puis, par une récurrence finie immédiate à partir des mêmes opérations sur lignes et colonnes :

$$\chi_{M_p} = (X-\beta_p)^{p-1} \left| X + (p-1)\beta_p \right|_{[p-(p-1)]} = (X-\beta_p)^{p-1}\big(X + (p-1)\beta_p\big)$$

On en déduit la diagonalisée $D_p = \mathrm{Diag}\Big((1-p)\beta_p\,;\beta_p\,;\dots;\beta_p\Big)$.

Or, puisque $\beta_p = {\alpha_p}^{p-2}$ et $0 < \alpha_p < 1$, on a : $\forall\, p \ge 3,\ 0 < \beta_p < 1$. Montrons qu'on a de plus : $\forall\, p \ge 3,\ (1-p)\beta_p < -1$. Soit $\theta_p$ la fonction définie sur $[0;1]$ par : $\theta_p(x) = x^{p-1} + x - 1$. Elle s'annule en $\alpha_p$ et on montre facilement qu'elle est strictement croissante. Or :

$$\forall\, p \geq 3,\ \ \theta_p\left(1-\frac{1}{p}\right) = \left(1-\frac{1}{p}\right)^{p-1} - \frac{1}{p} > 0$$

La démonstration du théorème 3 suivie dans (Pouvreau, Bouis, 2018) a en effet déjà établi que :

$$\forall\, p \geq 3,\ \ \frac{p-1}{p^{\frac{p-2}{p-1}}} > 1$$

Et cette inégalité est clairement équivalente à celle annoncée. On en déduit :

$$\forall\, p \geq 3,\ \ \alpha_p < 1 - \frac{1}{p}$$

Cette inégalité équivaut à son tour à :

$$\forall\, p \geq 3,\ \ (p-1)(1-\alpha_p) > \alpha_p$$

C'est-à-dire aussi à :

$$\forall\, p \geq 3,\ \ (p-1){\alpha_p}^{p-1} > \alpha_p$$

Et donc bien, en fin de compte, à :

$$\forall\, p \geq 3,\ \ (1-p)\beta_p < -1$$

Par conséquent, la suite extraite $\left(\left((1-p)\beta_p\right)^{2m}\right)_{m\in\mathbb{N}}$ diverge exponentiellement vers $+\infty$ et la suite extraite $\left(\left((1-p)\beta_p\right)^{2m+1}\right)_{m\in\mathbb{N}}$ diverge exponentiellement vers $-\infty$. Comme la suite $\left({\beta_p}^m\right)_{m\in\mathbb{N}}$ converge vers 0, les éléments de réduction précédents montrent que les suites extraites $\left({D_p}^{2m}\right)_{m\in\mathbb{N}}$ et $\left({D_p}^{2m+1}\right)_{m\in\mathbb{N}}$ divergent exponentiellement (dans des directions matricielles distinctes) au sens d'une norme quelconque sur $\mathcal{M}_p(\mathbb{R})$. On en déduit bien que $\Omega$ est exponentiellement instable.

## 2.2. Convergence exponentielle d'une suite duale vers le centre de gravité de $\mathcal{A}$

Nous allons maintenant établir dans toute sa généralité le théorème suivant :

**Théorème 2**
Soit $\mathfrak{B}$ une suite barypolygonale de $\mathcal{A}$, avec $\mathrm{Card}(\mathcal{A}) = p \geq 3$.
La suite duale de $\mathfrak{B}$ converge exponentiellement vers le centre de gravité de $\mathcal{A}$.

Démonstration.
Par souci de clarté, nous découperons la démonstration au moyen de lemmes intermédiaires. Le résultat étant acquis si $\mathfrak{B}$ est régulière, on supposera dans ce qui suit qu'elle ne l'est pas.

Comme dans le cas où $p = 3$ (où le théorème est déjà établi) et par le même argument, on peut supposer sans inconvénient que $0 < u_1^{(0)} \leq u_2^{(0)} \leq u_3^{(0)} \leq \cdots \leq u_p^{(0)} < 1$.

Soit $E$ une partie de l'ensemble $[\![1;p]\!]$. Dans tout ce qui suit, nous conviendrons de noter

$$Q_E^{(m)} = \prod_{\substack{1\leq r\leq p \\ r\notin E}} \left(X - u_r^{(m)}\right)$$

et, pour tout élément $i$ de l'ensemble $[\![0;p-\mathrm{Card}(E)]\!]$, de noter $\sigma_{i\,;E}^{(m)}$ la $i$-ème fonction symétrique élémentaire du polynôme $Q_E^{(m)}$.

**Lemme (a)**

$$\forall\, m \in \mathbb{N}\,,\ \ 0 < u_1^{(m)} \le u_2^{(m)} \le u_3^{(m)} \le \cdots \le u_p^{(m)} < 1$$

Vraie lorsque $m = 0$, cette propriété est héréditaire car pour tous $k \in [\![1; p-1]\!]$ et $m \in \mathbb{N}$ :

$$u_k^{(m)} - u_{k+1}^{(m)} = \left(1 - \prod_{\substack{1\le i\le p\\ i\ne k}} u_i^{(m)}\right) - \left(1 - \prod_{\substack{1\le i\le p\\ i\ne k+1}} u_i^{(m)}\right)$$

$$= \prod_{\substack{1\le i\le p\\ i\ne k+1}} u_i^{(m)} - \prod_{\substack{1\le i\le p\\ i\ne k}} u_i^{(m)} = \left(u_k^{(m)} - u_{k+1}^{(m)}\right) \prod_{\substack{1\le i\le p\\ i\ne k\,;\, i\ne k+1}} u_i^{(m)}$$

**Lemme (b)**

Pour tout $(k; l) \in [\![1; p]\!]^2$ tel que $k < l$, la suite $\left(u_l^{(2q)}/u_k^{(2q)}\right)_{q\in\mathbb{N}}$ est décroissante et minorée par 1, donc convergente vers un réel $L_{(k;l)} \ge 1$.

Commençons par observer que pour tout $k \in [\![1; p]\!]$ et tout $m \in \mathbb{N}$ :

$$u_k^{(m+2)} = 1 - \prod_{\substack{1\le i\le p\\ i\ne k}} u_i^{(m+1)} = 1 - \prod_{\substack{1\le i\le p\\ i\ne k}} \left(1 - \prod_{\substack{1\le j\le p\\ j\ne i}} u_j^{(m)}\right)$$

$$= u_k^{(m)} \sum_{0\le j\le p-2} (-1)^j {\pi_m}^j\, \sigma_{p-2-j\,;\,\{k\}}^{(m)}$$

où :

$$\pi_m = \prod_{1\le r\le p} u_r^{(m)}$$

Ce qui revient aussi à dire qu'on a pour tout $k \in [\![1; p]\!]$ :

$$u_k^{(m+2)} = \frac{1}{\pi_m} u_k^{(m)} \left(\sigma_{p-1\,;\,\{k\}}^{(m)} - Q_{\{k\}}^{(m)}(\pi_m)\right)$$

On obtient alors pour tout $(k; l) \in [\![1; p]\!]^2$ :

$$u_k^{(m)} u_l^{(m+2)} - u_l^{(m)} u_k^{(m+2)}$$

$$= u_k^{(m)} u_l^{(m)} \sum_{0\le j\le p-2} (-1)^j\, {\pi_m}^j \left(\sigma_{p-2-j\,;\,\{l\}}^{(m)} - \sigma_{p-2-j\,;\,\{k\}}^{(m)}\right)$$

$$= u_k^{(m)} u_l^{(m)} \left(u_k^{(m)} - u_l^{(m)}\right) \sum_{0\le j\le p-3} (-1)^j\, {\pi_m}^j\, \sigma_{p-3-j\,;\,\{k;l\}}^{(m)}$$

Il reste alors à observer qu'en posant $\gamma_{\{k;l\}}^{(m)} = u_k^{(m)} u_l^{(m)}$ , on a :

$$\sum_{0\le j\le p-3} (-1)^j\, {\pi_m}^j\, \sigma_{p-3-j\,;\,\{k;l\}}^{(m)} = \frac{1}{\gamma_{\{k;l\}}^{(m)}} \left(1 - \prod_{\substack{1\le r\le p\\ r\notin\{k;l\}}} \left(1 - \gamma_{\{k;l\}}^{(m)} u_r^{(m)}\right)\right)$$

Ce qui montre, puisque $0 < \gamma_{\{k;l\}}^{(m)} u_r^{(m)} < 1$ pour tout $r \in [\![1;p]\!] \backslash \{k;l\}$, que :

$$\sum_{0 \leq j \leq p-3} (-1)^j {\pi_m}^j \ \sigma_{p-3-j\,;\,\{k;l\}}^{(m)} > 0$$

On obtient ainsi, compte tenu du lemme (a) :

$$\forall\, (k;l) \in [\![1;p]\!]^2,\ (k<l) \Rightarrow \left(\forall\, m \in \mathbb{N},\ u_k^{(m)} u_l^{(m+2)} - u_l^{(m)} u_k^{(m+2)} \leq 0\right)$$

D'où aussi :

$$\forall\, (k;l) \in [\![1;p]\!]^2,\ (k<l) \Rightarrow \left(\forall\, m \in \mathbb{N},\ 1 \leq \frac{u_l^{(m+2)}}{u_k^{(m+2)}} \leq \frac{u_l^{(m)}}{u_k^{(m)}}\right)$$

L'assertion du lemme (b) est ainsi établie.

**Lemme (c)**

$$\forall\, m \in \mathbb{N},\ \frac{\dfrac{u_p^{(m+2)}}{u_1^{(m+2)}} - 1}{\dfrac{u_p^{(m)}}{u_1^{(m)}} - 1} < \frac{1}{2}$$

En reprenant les notations introduites précédemment, on pose d'abord, pour tout $m \in \mathbb{N}$ :

$$K_m = \prod_{1<r<p} \left(u_r^{(m)} - \pi_m\right) = (-1)^p Q_{\{1;p\}}^{(m)}(\pi_m)$$

et

$$C_m = u_1^{(m)} u_p^{(m)} \sum_{0 \leq j \leq p-3} (-1)^j \,{\pi_m}^j\, \sigma_{p-3-j\,;\,\{1;p\}}^{(m)} = u_1^{(m)} u_p^{(m)} \frac{\sigma_{p-2\,;\,\{1;p\}}^{(m)} + (-1)^{p-1} Q_{\{1;p\}}^{(m)}(\pi_m)}{\pi_m}$$

Remarquons que $K_m > 0$ pour tout $m \in \mathbb{N}$ en raison des définitions de $Q_{\{1;p\}}^{(m)}$ de $\pi_m$, et du fait que $u_k^{(m)} \in ]0;1[$ pour tout $(k;m) \in [\![1;p]\!] \times \mathbb{N}$ . Remarquons aussi que $C_m > 0$. En effet :

$$\sigma_{p-2\,;\,\{1;p\}}^{(m)} + (-1)^{p-1} Q_{\{1;p\}}^{(m)}(\pi_m) = \prod_{2 \leq j \leq p-1} u_j^{(m)} - \prod_{2 \leq j \leq p-1} (u_j^{(m)} - \pi_m)$$

avec $u_j^{(m)} > (u_j^{(m)} - \pi_m)$ pour tout $(j;m) \in [\![1;p]\!] \times \mathbb{N}$. Observons alors que :

$$\forall\, m \in \mathbb{N},\ \begin{cases} u_1^{(m+2)} = u_1^{(m)} \displaystyle\sum_{0 \leq j \leq p-2} (-1)^j {\pi_m}^j \ \sigma_{p-2-j\,;\,\{1\}}^{(m)} = K_m u_1^{(m)} + C_m \\ u_p^{(m+2)} = u_p^{(m)} \displaystyle\sum_{0 \leq j \leq p-2} (-1)^j {\pi_m}^j \ \sigma_{p-2-j\,;\,\{p\}}^{(m)} = K_m u_p^{(m)} + C_m \end{cases}$$

Ces deux relations permettent d'en déduire, l'irrégularité supposée de la suite barypolygonale initiale impliquant bien l'existence du premier rapport :

$$\forall\, m \in \mathbb{N},\ \frac{\dfrac{u_p^{(m+2)}}{u_1^{(m+2)}} - 1}{\dfrac{u_p^{(m)}}{u_1^{(m)}} - 1} = \frac{\dfrac{K_m u_p^{(m)} + C_m}{K_m u_1^{(m)} + C_m} - 1}{\dfrac{u_p^{(m)}}{u_1^{(m)}} - 1} = \frac{K_m u_1^{(m)}}{K_m u_1^{(m)} + C_m} = 1 - \frac{C_m}{K_m u_1^{(m)} + C_m}$$

Or, pour tout $m \in \mathbb{N}$ :

$$\frac{K_m u_1^{(m)} + C_m}{C_m} = 1 + \frac{K_m u_1^{(m)}}{C_m}$$

avec :

$$\frac{K_m u_1^{(m)}}{C_m} = \frac{(-1)^p Q_{\{1;p\}}^{(m)}(\pi_m)}{u_1^{(m)} u_p^{(m)} \frac{\sigma_{p-2\,;\,\{1;p\}}^{(m)} + (-1)^{p-1} Q_{\{1;p\}}^{(m)}(\pi_m)}{\pi_m}} u_1^{(m)}$$

$$= \frac{1}{u_p^{(m)}} \frac{(-1)^p Q_{\{1;p\}}^{(m)}(\pi_m)}{\sigma_{p-2\,;\,\{1;p\}}^{(m)} - (-1)^p Q_{\{1;p\}}^{(m)}(\pi_m)}$$

$$= u_1^{(m)} \frac{(-1)^p Q_{\{1;p\}}^{(m)}(\pi_m)\, \sigma_{p-2\,;\,\{1;p\}}^{(m)}}{\sigma_{p-2\,;\,\{1;p\}}^{(m)} - (-1)^p Q_{\{1;p\}}^{(m)}(\pi_m)}$$

$$= u_1^{(m)} \frac{\prod_{2 \le j \le p-1} u_j^{(m)} \times \prod_{2 \le j \le p-1} \left(u_j^{(m)} - \pi_m\right)}{\prod_{2 \le j \le p-1} u_j^{(m)} - \prod_{2 \le j \le p-1} \left(u_j^{(m)} - \pi_m\right)}$$

$$= u_1^{(m)} \frac{\prod_{2 \le j \le p-1} \left(u_j^{(m)}\right)^2 \times \prod_{2 \le j \le p-1} \left(1 - \prod_{1 \le i \le p\,;\, i \ne j} u_i^{(m)}\right)}{\prod_{2 \le j \le p-1} u_j^{(m)} - \prod_{2 \le j \le p-1} u_j^{(m)} \times \prod_{2 \le j \le p-1} \left(1 - \prod_{1 \le i \le p\,;\, i \ne j} u_i^{(m)}\right)}$$

$$= \prod_{1 \le j \le p-1} u_j^{(m)} \times \frac{\prod_{2 \le j \le p-1} \left(1 - \prod_{1 \le i \le p\,;\, i \ne j} u_i^{(m)}\right)}{1 - \prod_{2 \le j \le p-1} \left(1 - \prod_{1 \le i \le p\,;\, i \ne j} u_i^{(m)}\right)}$$

Il est donc déjà clair que :

$$\forall\, m \in \mathbb{N},\ 0 < \frac{K_m u_1^{(m)}}{C_m} < \prod_{1 \le j \le p-1} u_j^{(m)} \times \frac{1}{1 - \prod_{2 \le j \le p-1} \left(1 - \prod_{1 \le i \le p\,;\, i \ne j} u_i^{(m)}\right)}$$

De plus, on a aussi lorsque $p \geq 4$ :

$$\prod_{1 \le j \le p-1} u_j^{(m)} \times \frac{1}{1 - \prod_{2 \le j \le p-1} \left(1 - \prod_{1 \le i \le p\,;\, i \ne j} u_i^{(m)}\right)}$$

$$= \prod_{1 \le j \le p-1} u_j^{(m)} \times \frac{1}{\sum_{1 \le j \le p-2} (-1)^{j-1} \left(\gamma_{\{1;p\}}^{(m)}\right)^j \sigma_{j\,;\,\{1;p\}}^{(m)}}$$

$$= \frac{1}{u_p^{(m)}} \prod_{2 \le j \le p-1} u_j^{(m)} \times \frac{1}{\sum_{1 \le j \le p-2} (-1)^{j-1} \left(\gamma_{\{1;p\}}^{(m)}\right)^{j-1} \sigma_{j\,;\,\{1;p\}}^{(m)}}$$

Par conséquent, pour $p \geq 4$ :

$$\forall\, m \in \mathbb{N},\ 0 < \frac{K_m u_1^{(m)}}{C_m} < \frac{u_3^{(m)}}{u_p^{(m)}} \times \frac{1}{\frac{1}{u_2^{(m)}} \sum_{1 \le j \le p-2} (-1)^{j-1} \left(\gamma_{\{1;p\}}^{(m)}\right)^{j-1} \sigma_{j\,;\,\{1;p\}}^{(m)}}$$

Or, pour tout $m \in \mathbb{N}$ :

$$\frac{1}{u_2^{(m)}} \sum_{1\leq j\leq p-2} (-1)^{j-1} \left(\gamma_{\{1;p\}}^{(m)}\right)^{j-1} \sigma_{j\,;\{1;p\}}^{(m)} = 1 + \left(\frac{1}{u_2^{(m)}} - \gamma_{\{1;p\}}^{(m)}\right) \beta_{\{1;p\}}^{(m)}$$

avec :

$$\beta_{\{1;p\}}^{(m)} = \sum_{1\leq j\leq p-3} (-1)^{j-1} \left(\gamma_{\{1;p\}}^{(m)}\right)^{j-1} \sigma_{j\,;\{1;2;p\}}^{(m)}$$

Remarquons alors que :

$$\beta_{\{1;p\}}^{(m)} = \frac{1}{\gamma_{\{1;p\}}^{(m)}} - \left(\gamma_{\{1;p\}}^{(m)}\right)^{p-4} Q_{\{1;2;p\}}^{(m)}\left(\frac{1}{\gamma_{\{1;p\}}^{(m)}}\right) = \frac{1}{\gamma_{\{1;p\}}^{(m)}}\left(1 - \left(\gamma_{\{1;p\}}^{(m)}\right)^{p-3} Q_{\{1;2;p\}}^{(m)}\left(\frac{1}{\gamma_{\{1;p\}}^{(m)}}\right)\right)$$

Il reste finalement à observer que :

$$\left(\gamma_{\{1;p\}}^{(m)}\right)^{p-3} Q_{\{1;2;p\}}^{(m)}\left(\frac{1}{\gamma_{\{1;p\}}^{(m)}}\right) = \left(\gamma_{\{1;p\}}^{(m)}\right)^{p-3} \prod_{3\leq j\leq p-1}\left(\frac{1}{\gamma_{\{1;p\}}^{(m)}} - u_j^{(m)}\right) = \prod_{3\leq j\leq p-1}\left(1 - \gamma_{\{1;p\}}^{(m)}\, u_j^{(m)}\right)$$

Ce qui implique :

$$\left(\gamma_{\{1;p\}}^{(m)}\right)^{p-3} Q_{\{1;2;p\}}^{(m)}\left(\frac{1}{\gamma_{\{1;p\}}^{(m)}}\right) < 1$$

On peut en déduire : $\forall\, m \in \mathbb{N},\ \beta_{\{1;p\}}^{(m)} > 0$. Et donc aussi, puisque

$$\left(\forall\, m \in \mathbb{N},\ \frac{1}{u_2^{(m)}} > 1 > \gamma_{\{1;p\}}^{(m)}\right) \Rightarrow \left(\forall\, m \in \mathbb{N},\ \frac{1}{u_2^{(m)}} - \gamma_{\{1;p\}}^{(m)} > 0\right)$$

que :

$$\forall\, m \in \mathbb{N},\ \frac{1}{u_2^{(m)}} \sum_{1\leq j\leq p-2} (-1)^{j-1} \left(\gamma_{\{1;p\}}^{(m)}\right)^{j-1} \sigma_{j\,;\{1;p\}}^{(m)} > 1$$

On obtient donc finalement :

$$\forall\, m \in \mathbb{N},\ ,\ 0 < \frac{K_m u_1^{(m)}}{C_m} < \frac{u_3^{(m)}}{u_p^{(m)}}$$

Il résulte alors du lemme (a) que :

$$\forall\, m \in \mathbb{N},\quad 0 < \frac{K_m u_1^{(m)}}{C_m} < 1$$

D'où :

$$\forall\, m \in \mathbb{N},\ \frac{C_m}{K_m u_1^{(m)} + C_m} > \frac{1}{2}$$

Puis en définitive les inégalités annoncées.

**Lemme (d)**

$$\forall\,(k;l) \in [\![1;p]\!]^2,\ \lim_{m\to+\infty} \frac{u_l^{(m)}}{u_k^{(m)}} = 1$$

(cette convergence étant de type exponentiel)

Il résulte en effet directement par récurrence du lemme (c) que :

$$\forall\, q \in \mathbb{N},\ \ 0 < \frac{u_p^{(2q)}}{u_1^{(2q)}} - 1 < \left(\frac{1}{2}\right)^q \left(\frac{u_p^{(0)}}{u_1^{(0)}} - 1\right)$$

D'où $L_{(1;p)} = 1$. On en déduit ensuite par encadrement que pour tout $(k; l) \in [\![1; p]\!]^2$ tel que $k < l$, $L_{(k;l)} = 1$, les convergences étant là aussi de type exponentiel. On montre de même que la suite $\left(u_l^{(2q+1)}/u_k^{(2q+1)}\right)_{q\in\mathbb{N}}$ converge exponentiellement vers 1 pour tout $(k; l) \in [\![1; p]\!]^2$. D'où le résultat.

**Terme de la démonstration**

Observons que :

$$\forall\, k \in [\![1; p]\!],\ \ \forall\, q \in \mathbb{N},\ \ t_k^{(2q+1)} = \prod_{\substack{1\le i\le p \\ i\ne k}} \left(1 - t_i^{(2q)}\right) = \prod_{\substack{1\le i\le p \\ i\ne k}} u_i^{(2q)}$$

Cette égalité implique :

$$\forall\, (k; l) \in [\![1; p]\!]^2,\ \ \forall\, q \in \mathbb{N},\qquad \frac{t_l^{(2q+1)}}{t_k^{(2q+1)}} = \frac{u_k^{(2q)}}{u_l^{(2q)}}$$

Et donc d'après le lemme (d) :

$$\lim_{q\to+\infty} \frac{t_l^{(2q+1)}}{t_k^{(2q+1)}} = 1$$

La convergence étant de type exponentiel. On montre bien sûr de même que :

$$\lim_{q\to+\infty} \frac{t_l^{(2q)}}{t_k^{(2q)}} = 1$$

Avec là aussi une convergence de type exponentiel. Il en résulte finalement :

$$\forall\, (k; l) \in [\![1; p]\!]^2,\ \ \lim_{m\to+\infty} \frac{t_l^{(m)}}{t_k^{(m)}} = 1$$

Avec dans tous les cas une convergence de type exponentiel. Ceci achève d'établir le théorème 2.

Les figures 1 à 4 illustrent successivement, en dimension 2 et dans un cas où $p = 5$, les dérivées d'ordre 0 à 3 d'une suite $(0{,}3; 0{,}08; 0{,}06; 0{,}04; 0{,}01)$-barypolygonale d'un ensemble $\mathcal{A}$ :

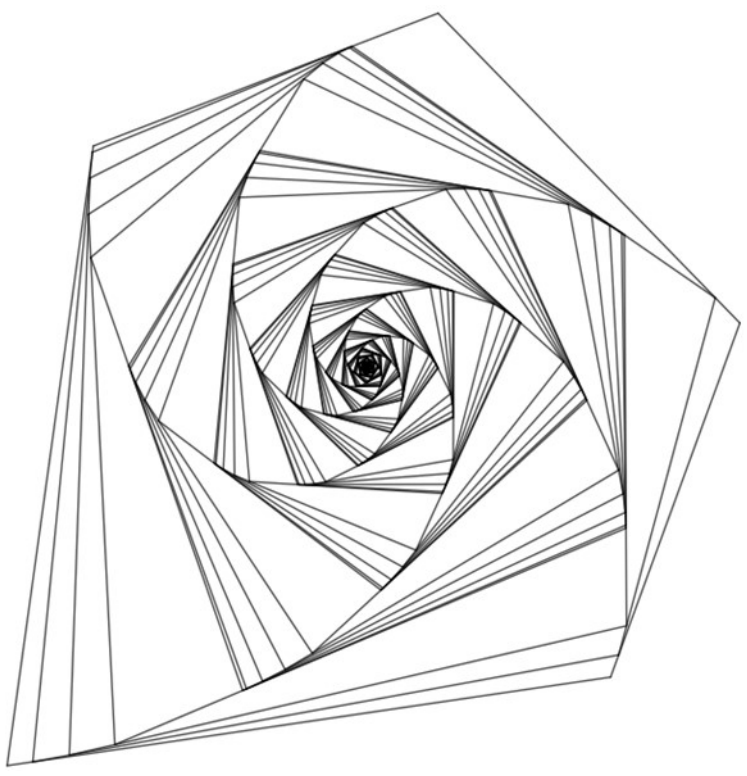

Figure 1

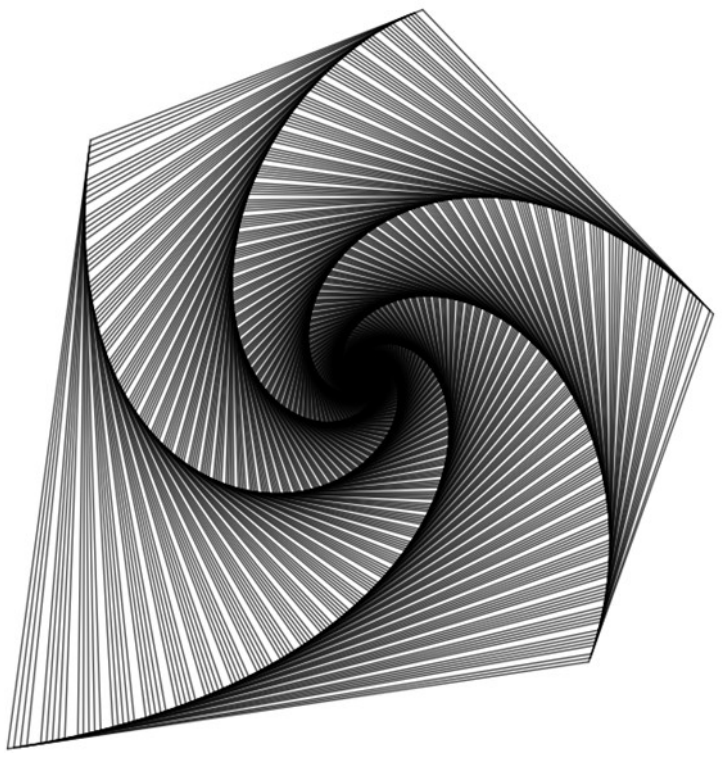

Figure 2

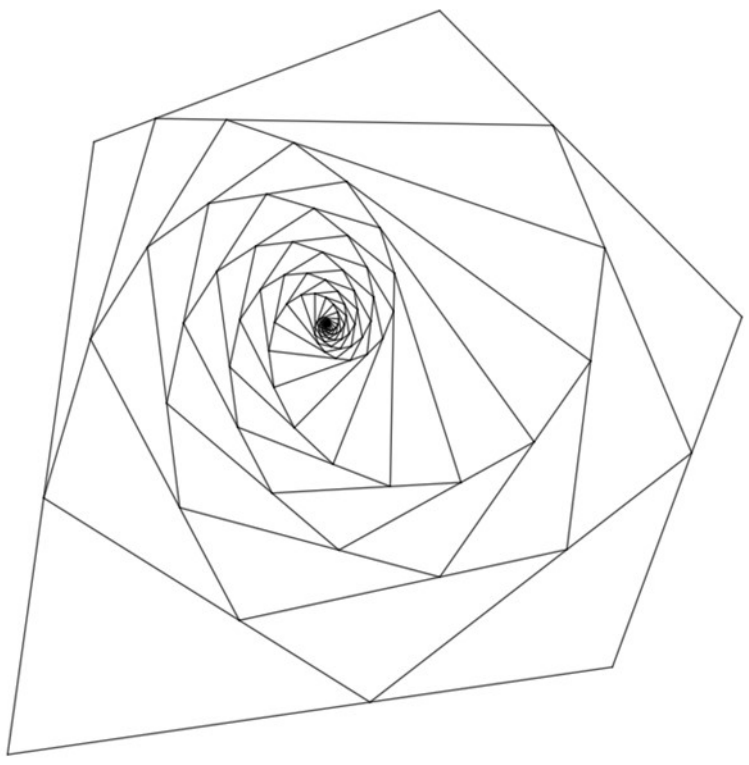

Figure 3

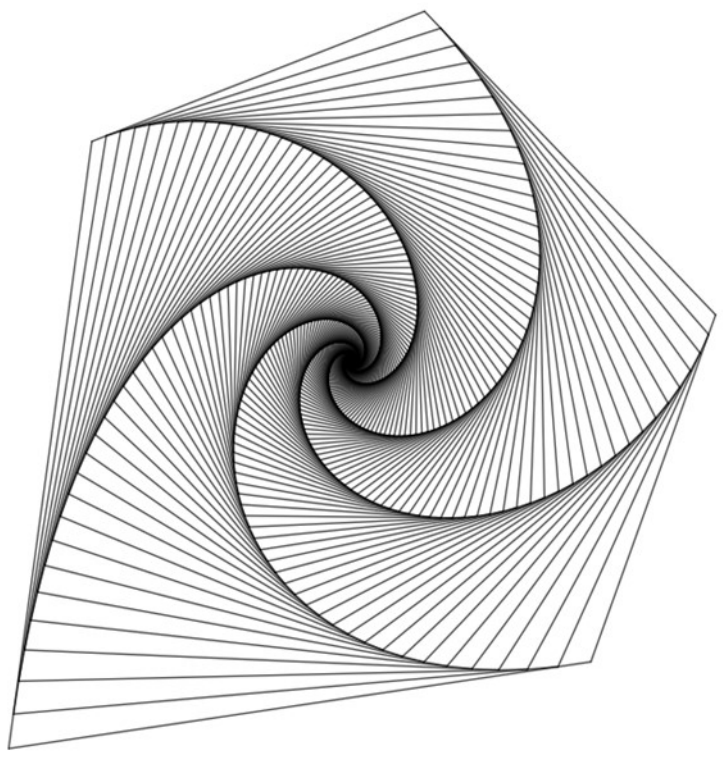

Figure 4

La superposition sur une même figure d'une succession de suites barypolygonales dérivées donne à voir une suite d'attracteurs convergeant vers le centre de gravité de $\mathcal{A}$, mais la convergence de type exponentiel rend difficile sa réalisation de manière lisible. La figure 5 fournit une illustration avec la superposition des deux suites dérivées d'ordre 0 et 1 lorsque $t = (0{,}03; 0{,}02; 0{,}03; 0{,}02; 0{,}01)$.

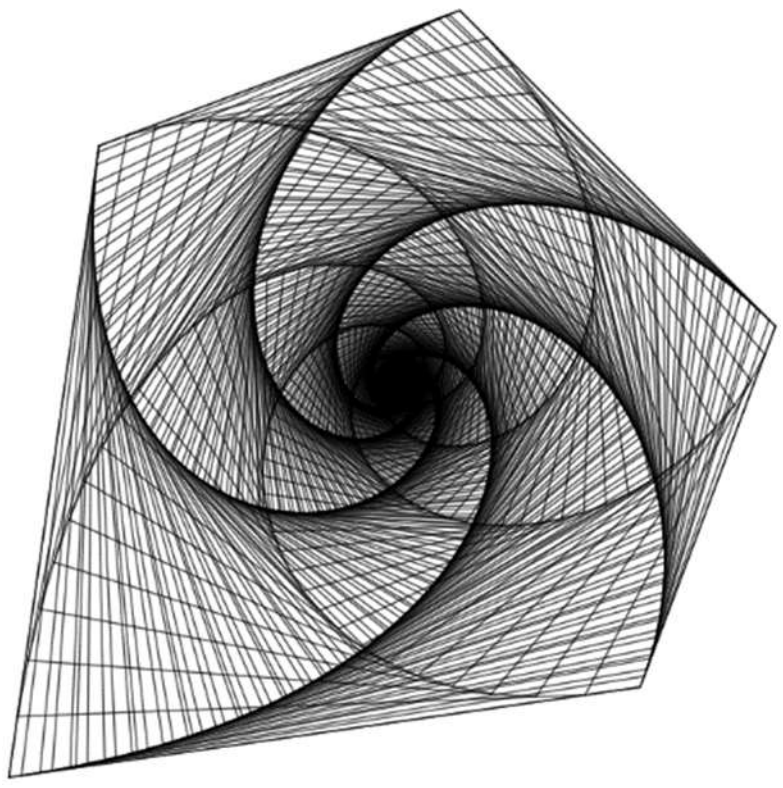

Figure 5

## 2.3. Lieu des éléments du système barypolygonal dérivé

On peut désormais généraliser l'étude réalisée lorsque $p = 3$ concernant le lieu des termes d'un $p$-uplet $\left(u_k^{(m)}\right)_{1\le k\le p}$ par rapport au point stationnaire « répulsif » $\left(\alpha_p\right)_{1\le k\le p}$ :

**Théorème 3**

Soit $\mathfrak{B}$ une suite $t$-barypolygonale de $\mathcal{A}$, avec $\mathrm{Card}(\mathcal{A}) = p \ge 3$. Soit $\alpha_p$ l'unique solution dans $[0;1]$ de l'équation $\left(E_p\right)$ : $x^{p-1} + x - 1 = 0$. Il existe $m_0 \in \mathbb{N}$ tel que pour tout $q \in \mathbb{N}$ :

$$\left(t_k^{(m_0+2q)}\right)_{1\le k\le p} \in \left]0; 1-\alpha_p\right[^p \text{ et } \left(t_k^{(m_0+2q+1)}\right)_{1\le k\le p} \in \left]1-\alpha_p; 1\right[^p$$

ou

$$\left(t_k^{(m_0+2q)}\right)_{1\le k\le p} \in \left]1-\alpha_p; 1\right[^p \text{ et } \left(t_k^{(m_0+2q+1)}\right)_{1\le k\le p} \in \left]0; 1-\alpha_p\right[^p$$

Démonstration.

On reprend ici l'hypothèse sans inconvénient selon laquelle

$$0 < u_1^{(0)} \le u_2^{(0)} \le u_3^{(0)} \le \cdots \le u_p^{(0)} < 1$$

On a alors vu que :

$$\forall\, m \in \mathbb{N}\,,\ \ 0 < u_1^{(m)} \leq u_2^{(m)} \leq u_3^{(m)} \leq \cdots \leq u_p^{(m)} < 1$$

D'après le théorème de Bolzano-Weierstrass, la suite $\left(\left(u_k^{(m)}\right)_{1\leq k\leq p}\right)_{m\in\mathbb{N}}$ d'éléments du compact $[0;1]^p$ admet une suite extraite $\left(\left(u_k^{(\psi(m))}\right)_{1\leq k\leq p}\right)_{m\in\mathbb{N}}$ convergeant vers un $p$-uplet $(\ell_k)_{1\leq k\leq p}$ de $[0;1]^p$, qui est ici tel que $0 \leq \ell_1 \leq \cdots \leq \ell_p \leq 1$. Comme on a supposé le $p$-uplet initial $t$ distinct du $p$-uplet $\left(1-\alpha_p\right)_{1\leq k\leq p}$ et puisque le point stationnaire $\left(\alpha_p\right)_{1\leq k\leq p}$ de $(\Sigma)$ est exponentiellement instable (théorème 1), on a $(\ell_k)_{1\leq k\leq p} \neq \left(\alpha_p\right)_{1\leq k\leq p}$. Supposons que $\ell_1 < \alpha_p$. Il résulte du théorème 2 que :

$$\forall\, k \in [\![1;p-1]\!],\ \ \lim_{m\to+\infty} \frac{u_k^{(\psi(m))}}{u_{k+1}^{(\psi(m))}} = 1$$

On en déduit que si $\ell_1 = 0$, alors $\ell_k = 0$ pour tout $k \in [\![2;p]\!]$. Et que si $\ell_1 > 0$, alors on a encore $\ell_k = \ell_1$ pour tout $k \in [\![1;p]\!]$. Dans les deux cas, on peut donc écrire :

$$\forall\, \varepsilon > 0,\ \exists\, M \in \mathbb{N},\ \forall\, m \geq M,\ \ 0 < u_1^{(\psi(m))} \leq u_2^{(\psi(m))} \leq \cdots \leq u_p^{(\psi(m))} < \ell_1 + \varepsilon$$

En choisissant $\varepsilon = \frac{\alpha_p - \ell_1}{2}$ , on obtient ainsi :

$$\forall\, m \geq M,\ \ 0 < u_1^{(\psi(m))} \leq u_2^{(\psi(m))} \leq \cdots \leq u_p^{(\psi(m))} < \frac{\ell_1 + \alpha_p}{2} < \alpha_p$$

L'existence d'un entier $m_0 = \psi(M)$ satisfaisant la condition recherchée est de la sorte assurée.

On montre de même l'existence d'un tel entier lorsqu'on suppose $\ell_1 \geq \alpha_p$, $m_0$ étant cette fois tel que $\left(u_k^{(m_0)}\right)_{1\leq k\leq p} \in \left]\alpha_p;1\right[^p$ : il suffit alors d'utiliser le fait que nécessairement, $\alpha_p < \ell_p \leq 1$ et d'en déduire des inégalités analogues à celles qui précèdent sur l'intervalle $\left]\alpha_p;1\right[$.

On peut dès lors en déduire le résultat annoncé dans le théorème 3. Supposons par exemple que $m_0 \in \mathbb{N}$ soit tel que $\left(u_k^{(m_0)}\right)_{1\leq k\leq p} \in \left]\alpha_p;1\right[^p$. On montre alors par récurrence (avec ici encore deux formes de justification) que pour tout $k \in [\![1;p]\!]$ fixé et tout $q \in \mathbb{N}$ :

$$\left(u_k^{(m_0+2q)}\right)_{1\leq k\leq p} \in \left]\alpha_p;1\right[^p \ \text{ et } \ \left(u_k^{(m_0+2q+1)}\right)_{1\leq k\leq p} \in \left]0;\alpha_p\right[^p$$

Cette propriété est initialisée. Compte tenu de $1 - {\alpha_p}^{p-1} = \alpha_p$, on a en effet :

$$u_k^{(m_0+1)} = 1 - \prod_{\substack{1\leq j\leq p \\ j\neq k}} u_j^{(m_0)} = 1 - \prod_{\substack{1\leq j\leq p \\ j\neq k}} \left(\left(u_j^{(m_0)} - \alpha_p\right) + \alpha_p\right)$$

$$= \alpha_p - \sum_{1\leq j\leq p-1} s_{j\,;\{k\}}^{(m)}\, {\alpha_p}^{p-1-j} < \alpha_p$$

où $s_{j\,;\{k\}}^{(m)}$ est, pour $(k;m;j) \in [\![1;p]\!] \times \mathbb{N} \times [\![0;p-1]\!]$, la $j$-ème fonction symétrique élémentaire de

$$P_{\{k\}}^{(m)} = Q_{\{k\}}^{(m)}\left(X + \alpha_p\right) = \prod_{\substack{1\leq r\leq p \\ r\neq k}} \left(X - \left(u_r^{(m)} - \alpha_p\right)\right)$$

Il est alors clair, compte tenu de $\left(u_k^{(m_0)}\right)_{1\le k\le p} \in \left]\alpha_p;1\right[^p$ et de $\alpha_p > 0$, que $u_k^{(m_0+1)} < \alpha_p$ pour tout $k \in [\![1;p]\!]$. Soit maintenant $q \in \mathbb{N}$ tel que $\left(u_k^{(m_0+2q)}\right)_{1\le k\le p} \in \left]\alpha_p;1\right[^p$. On montre alors comme précédemment que $\left(u_k^{(m_0+2q+1)}\right)_{1\le k\le p} \in \left]0;\alpha_p\right[^p$. De plus, pour tout $k \in [\![1;p]\!]$ :

$$u_k^{(m_0+2q+2)} = 1 - \prod_{\substack{1\le j\le p\\ j\ne k}} u_j^{(m_0+2q+1)}$$

Or, pour tout $j \in [\![1;p]\!]\backslash\{k\}$, $0 < u_j^{(m_0+2q+1)} < \alpha_p$ ; ce qui implique :

$$0 < \prod_{\substack{1\le j\le p\\ j\ne k}} u_j^{(m_0+2q+1)} < {\alpha_p}^{p-1}$$

et donc :

$$\alpha_p = 1 - {\alpha_p}^{p-1} < 1 - \prod_{\substack{1\le j\le p\\ j\ne k}} u_j^{(m_0+2q+1)} < 1$$

Par conséquent : $u_k^{(m_0+2q+2)} \in \left]\alpha_p;1\right[$. On peut alors de nouveau établir comme dans le cas où $q = 0$ que $\left(u_k^{(m_0+2q+3)}\right)_{1\le k\le p} \in \left]0;\alpha_p\right[^p$. Ceci achève d'établir par récurrence la propriété annoncée.

De même, s'il existe $m_0 \in \mathbb{N}$ tel que $\left(u_k^{(m_0)}\right)_{1\le k\le p} \in \left]0;\alpha_p\right[^p$, alors pour tout $q \in \mathbb{N}$ :

$$\left(u_k^{(m_0+2q)}\right)_{1\le k\le p} \in \left]0;\alpha_p\right[^p \text{ et } \left(u_k^{(m_0+2q+1)}\right)_{1\le k\le p} \in \left]\alpha_p;1\right[^p$$

Toutes ces conclusions se transfèrent par « conjugaison » aux éléments du système dérivé $(S)$.

## 2.4. Convergence des suites extraites des termes de rang pair et de rang impair d'éléments du système dérivé

On peut désormais achever cette étude générale en démontrant la convergence des suites extraites des termes de rang pair et de rang impair d'éléments du système dérivé :

**Théorème 4**
Soit $\mathfrak{B}$ une suite $t$-barypolygonale de $\mathcal{A}$, avec $\mathrm{Card}(\mathcal{A}) = p \ge 3$. Les deux suites extraites $\left(\left(t_k^{(2m)}\right)_{1\le k\le p}\right)_{m\in\mathbb{N}}$ et $\left(\left(t_k^{(2m+1)}\right)_{1\le k\le p}\right)_{m\in\mathbb{N}}$ de $]0;1[^p$ convergent, l'une convergeant vers $(0)_{1\le k\le p}$ et l'autre vers $(1)_{1\le k\le p}$.

Démonstration.
L'étude réalisée au 2.3. montre qu'il suffit, sans perte de généralité, d'examiner le cas d'un entier $m_0 \in \mathbb{N}$ tel que $0 < u_1^{(m_0)} \le u_2^{(m_0)} \le \cdots \le u_p^{(m_0)} < \alpha_p$ ; ce cas impliquant que :

$$\forall\, m \ge m_0\,,\ \ 0 < u_1^{(m)} \le u_2^{(m)} \le \cdots \le u_p^{(m)} < \alpha_p$$

Montrons alors que les $p$ suites $\left(u_k^{(m_0+2q)}\right)_{q\in\mathbb{N}}$ convergent vers 0, cependant que les $p$ suites $\left(u_k^{(m_0+2q+1)}\right)_{q\in\mathbb{N}}$ convergent vers 1.

Considérons à cette fin la fonction $f_p : x \mapsto 1 - x^{p-1}$. C'est une bijection décroissante de $[0;1]$ dans $[0;1]$. Si $u_1^{(m_0+1)} > f_p(u_p^{(m_0)})$, on pose $\tau = u_p^{(m_0)}$ ; sinon, on pose $\tau = f_p{}^{-1}(u_1^{(m_0+1)})$. On définit alors la suite $(\tau_m)_{m \geq m_0}$ par :

$$\begin{cases} \tau_{m_0} = \tau \\ \forall\, m \geq m_0\,,\ \ \tau_{m+1} = f_p(\tau_m) \end{cases}$$

Montrons par récurrence que :

$$\text{Pour tout } q \in \mathbb{N},\ \ \tau_{m_0+2q} \geq u_p^{(m_0+2q)} \ \text{ et } \ \tau_{m_0+2q+1} \leq u_1^{(m_0+2q+1)}.$$

Ces deux inégalités sont en effet vraies pour $q = 0$, par définition de $\tau$. La première l'est évidemment. Quant à la seconde, elle découle du fait que

$$\left(u_1^{(m_0+1)} > f_p\left(u_p^{(m_0)}\right)\right) \Rightarrow \left(u_1^{(m_0+1)} > f_p(\tau) = f_p(\tau_{m_0}) = \tau_{m_0+1}\right)$$

$$\text{et } \left(u_1^{(m_0+1)} \leq f_p\left(u_p^{(m_0)}\right)\right) \Rightarrow \left(u_1^{(m_0+1)} = f_p(\tau) = \tau_{m_0+1}\right)$$

De plus, si on considère un entier $q \in \mathbb{N}$ tel que $\tau_{m_0+2q} \geq u_p^{(m_0+2q)}$ , alors :

$$u_1^{(m_0+2q+1)} = f_p\left(\sqrt[(p-1)]{\prod_{2\leq j\leq p} u_j^{(m_0+2q)}}\right) \geq f_p\left(u_p^{(m_0+2q)}\right) \geq f_p\left(\tau_{m_0+2q}\right) = \tau_{m_0+2q+1}$$

et on obtient de manière similaire que $\tau_{m_0+2q+2} \geq u_p^{(m_0+2q+2)}$, puis que $\tau_{m_0+2q+3} \leq u_1^{(m_0+2q+3)}$. Donc la propriété annoncée est aussi héréditaire.

Or, d'après le théorème 3 démontré dans (Pouvreau, Bouis, 2018), la suite $(\tau_{2m})$ converge vers 0 et la suite $(\tau_{2m+1})$ converge vers 1. Donc par comparaison, on déduit de ce qui précède que :

$$\lim_{m\to+\infty} u_1^{(2m+1)} = 1 \ \text{ et } \ \lim_{m\to+\infty} u_p^{(2m)} = 0$$

Ce qui implique, là encore par comparaison, qu'on a aussi :

$$\forall\, k \in [\![2;p]\!],\ \lim_{m\to+\infty} u_k^{(2m+1)} = 1 \quad \text{et} \quad \forall\, k \in [\![1;p-1]\!],\ \lim_{m\to+\infty} u_k^{(2m)} = 0$$

On obtient des résultats analogues si l'on suppose $m_0$ tel que $\alpha_p < u_1^{(m_0)} \leq u_2^{(m_0)} \leq \cdots \leq u_p^{(m_0)} < 1$.

Le comportement du système $(\Sigma)$ est ainsi déterminé dans tous les cas. Le résultat annoncé concernant le système $(S)$ s'en déduit alors immédiatement par « conjugaison ».

## 3. Conclusion

Cette étude a mis en évidence le caractère « régularisateur » de l'algorithme de dérivation d'une suite barypolygonale tel qu'il a été défini, en montrant que le centre de gravité du système de points initial concerné $\mathcal{A}$ de cardinal $p$ est un attracteur des suites dérivées lorsque $p \geq 3$.

On peut être tenté d'affirmer que par dérivation, toute suite barypolygonale irrégulière tend vers une suite barypolygonale régulière. Il est plus correct d'affirmer que par dérivation, toute suite barypolygonale irrégulière engendre deux suites extraites complémentaires de suites dérivées dont l'une converge vers la « pseudo-suite » régulière 0-barypolygonale de $\mathcal{A}$, et l'autre vers la « pseudo-suite » régulière 1-barypolygonale de $\mathcal{A}$ (ces deux « pseudo-suites » étant donc définies comme les deux valeurs d'adhérence de la suite barypolygonale dérivée). Les convergences dont il est ici question étant entendues au sens de la distance $d$ (ou de toute distance équivalente) définie sur l'ensemble $\mathfrak{B}_{\mathcal{A}}$ des suites barypolygonales de $\mathcal{A}$ complété par ces deux « pseudo-suites » régulières par :

$$\forall\, (\mathfrak{B}_1;\mathfrak{B}_2) \in \boldsymbol{\mathfrak{B}_{\mathcal{A}}}^2,\ d(\mathfrak{B}_1;\mathfrak{B}_2) = \mathrm{Sup}_{(k;n)\in[\![1;p]\!]\times\mathbb{N}}\left\{\left|t_{(1;k)}^{(n)} - t_{(2;k)}^{(n)}\right|\right\},$$

$\mathfrak{B}_1$ étant la suite $\left(t_{(1;1)}; t_{(1;2)}; \dots ; t_{(1;p)}\right)$-barypolygonale de $\mathcal{A}$ et $\mathfrak{B}_2$ la suite $\left(t_{(2;1)}; t_{(2;2)}; \dots ; t_{(2;p)}\right)$-barypolygonale de $\mathcal{A}$, $\left(t_{(1;k)}^{(n)}\right)_{n\in\mathbb{N}}$ et $\left(t_{(2;k)}^{(n)}\right)_{n\in\mathbb{N}}$ désignant ici les suites de coefficients définissant les termes successifs respectifs des suites $\mathfrak{B}_1$ et $\mathfrak{B}_2$ (notations utilisées dans (Pouvreau, Eupherte, 2016)).

# 4. Bibliographie